# Lamarckism and mechanism synthesis: approaching constrained optimization with ideas from biology

Wei Zhang[1,*], Xudong Shi[2], and Liwen Wang[2]

[1]Airport College, Civil Aviation University of China, Tianjin 300300, P.R. China.
[2]Aviation Ground Special Equipments Research Base, Civil Aviation University of China, Tianjin 300300, P.R. China.
[*]To whom correspondence should be addressed. E-mail: drwadecheung@gmail.com; telephone: +86-22-24092477

Abstract: Nonlinear constrained optimization problems are encountered in many scientific fields. To utilize the huge calculation power of current computers, many mathematic models are also rebuilt as optimization problems. Most of them have constrained conditions which need to be handled. Borrowing biological concepts, a study is accomplished for dealing with the constraints in the synthesis of a four-bar mechanism. Biologically regarding the constrained condition as a form of selection for characteristics of a population, four new algorithms are proposed, and a new explanation is given for the penalty method. Using these algorithms, three cases are tested in differential-evolution based programs. Better, or comparable, results show that the presented algorithms and methodology may become common means for constraint handling in optimization problems.

**I. Introduction**

Optimization is a process of finding the best solution to a problem from many alternatives. To ensure a fast and effective searching process, a skill or a set of steps called an optimization method is needed.

A recent trend in research on optimization methods is organizing the optimization process with ideas from other scientific fields, instead of only working in mathematics. Among these ideas, there have sorts of method called evolutionary algorithms (EAs) [1]. Holland [2] suggested that an optimization problem can be considered an evolutionary process. A population (composed of feasible solution vectors) gradually evolves towards the optimal solution by propagation, mutation, and selection. This is known as a "genetic algorithm" (GA). Methods that came from similar sources include: clonal algorithms [3], particle swarm optimization [4] (PSO), and differential evolution [5] (DE).

Most realistic optimization problems are multi-dimensional. Many of them are confined by conditions which are called constraints. Dealing with constraints is always a vital but difficult step in an optimization problem. Scientists have applied many algorithms for it, including penalty methods [6-11], adaptive constraints [12], co-evolution [13], multi-objective optimization [14], culture algorithms [15], Pareto method [16-17], and the augmented Lagrangian method [18-19]. The penalty method is the earliest constraint handling method, still receiving the most attention. Currently, no universal method has been developed.

It is surprising that no one has tried to explain handling constraints using the concepts of evolution. It will be simpler and more integrated if optimization can be organized as an evolution with self-contained constraint handling. Instead, we now use EAs (or other methods) together with other expertise for constraint handling to organize the optimization process.

One further step has to be taken for optimization to be described as an entirely evolutionary process: the constraints need to be expressed in biological terms by going back to evolution. New mechanisms to handle the constraints of nonlinear constrained optimization will be proposed. As an illustration, a classic optimization problem in engineering, the synthesis of a four-bar mechanism, will be investigated.

**II. Synthesis of the mechanism with optimization method**

**II.1 Synthesis of a mechanism**

A four-bar mechanism, a simple and classic mechanism, consists of four main bars linked with each other end to end. (See Figure 1)

Each bar can be denoted by a planar vector. Coordinates of the endpoint $E$ in a fixed reference system {O} can be expressed as following,

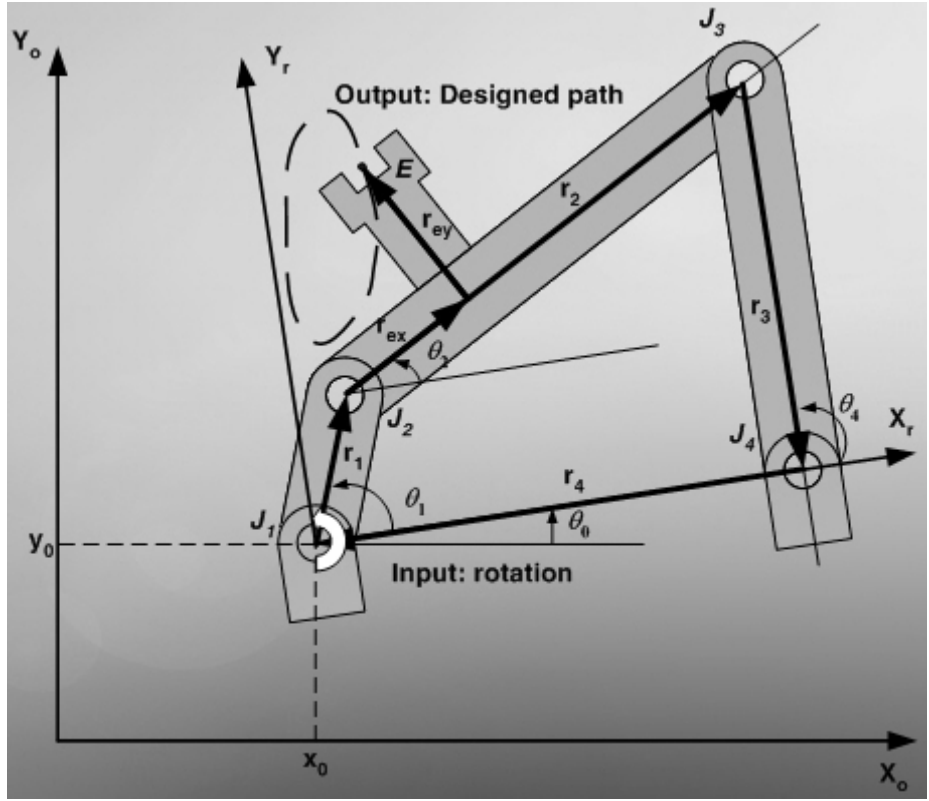

**Fig. 1.** Sketch of a four-bar mechanism. A desired motion path of the endpoint $E$ can be realized from an input rotation in joint $J_1$ via the mechanism. The fourth bar $r_4$ always represents fixed ground.

$$\left.\begin{aligned} X &= r_1 \cos(\theta_0 + \theta_1) + r_{ex} \cos(\theta_0 + \theta_2) - r_{ey} \sin(\theta_0 + \theta_2) + x_0 \\ Y &= r_1 \sin(\theta_0 + \theta_1) + r_{ex} \sin(\theta_0 + \theta_2) + r_{ey} \cos(\theta_0 + \theta_2) + y_0 \end{aligned}\right\} = f(\text{tht1}) \quad (1)$$

Considering the mechanism as a closed quadrangle, there is a group of additional equations,

$$\left.\begin{aligned} r_1 \cos\theta_1 + r_2 \cos\theta_2 - r_3 \cos\theta_4 - r_4 &= 0 \\ r_1 \sin\theta_1 + r_2 \sin\theta_2 - r_3 \sin\theta_4 &= 0 \end{aligned}\right\} \quad (2)$$

If the parameters of the mechanism $[r_1, r_2, r_3, r_4, r_{ex}, r_{ey}, \theta_0, x_0, y_0]$ are known, $\theta_2$ can be eliminated from equations (1) by using (2). Thus, the position of the endpoint is a function related to the input angle $\theta_1$. Solving for output $[X, Y]$ according to input angle $\theta_1$ is called forward kinematics. Otherwise, it is called inverse kinematics.

Synthesis of the mechanism is finding the parameters of mechanism for the desired output positions. It is hard to get one solution from so many possible solutions because there are too many variables in a limited number of equations. Sometimes, input angle $\theta_1$ is wanted.

Han [20] tried to solve this problem using numerical method. A local optimal solution is achieved through optimizing an objective function concerning distances between designed points and calculated points. Fang [21] first introduced GA into the synthesis of mechanisms. Optimization of mechanisms is treated as the evolution of a population.

Till (22) developed a distance function with the following format,

$$\sum_{i=1}^{N}\left[\left(X_d^i - X^i\right)^2 + \left(Y_d^i - Y^i\right)^2\right] \tag{3}$$

Where:

$N$ is the number of desired points;

$[X^i, Y^i]$ is the position of the $i$-th point, calculated with equations (1);

$[X_d^i, Y_d^i]$ is the desired position of $i$-th point.

Ideally, a solution of parameters can make equation (3) zero. It means that a mechanism with the parameters obtained leads to the endpoint passing through every desired point without error. An optimal solution is also acceptable, if (3) is made small within a prescribed range.

A suitable optimization method is needed to minimize the objective function (3).

A usable mechanism has to meet two basic constraint conditions. They are:

a. To ensure circulating motion, a four-bar mechanism should be capable of doing an entire rotation. This is known as the Grashof condition,

$$r_1 + r_4 < r_2 + r_3, \text{ while } r_4 < r_3 < r_2 < r_1 \tag{4}$$

A set of points on the designed path should be passed through in order.

**II.2 Differential evolution, as an optimization method**

Current methods for synthesis of mechanism are EAs, including: PSO, GA, and DE. Studies [16-17, 23-25] with DE are very successful. A DE method, proposed by Storn and Price [5], has already been applied to many fields [6, 26] due to its simplicity and superior performance compared with GA.

EAs are all built with reference to the natural selection of Darwinism. The optimal result will be achieved gradually from an evolved population by calculating fitness of every solution in every generation.

The optimization flow of DE method can be expressed as following,

**III. Constraints handling--- using biological concepts.**

Kunjur and Krishnamurty [22] reported that real number optimization performed better than binary optimization. They dealt with constraints through modifying the crossover operator and the mutation operator. Cabrera [25] achieved better results with a penalty method based on DE. Other methods which have been used for handling constraints include: penalty improvement [6, 9], Pareto improvement [14, 16-17], and customized factor [12, 22].

How can we understand constraint as a biological concept, to make the optimization process simpler to understand and more effective to execute? Constraints have not been taken into account in this way in the past.

In Figure 2, comparisons are drawn between the concepts of optimization and evolution.

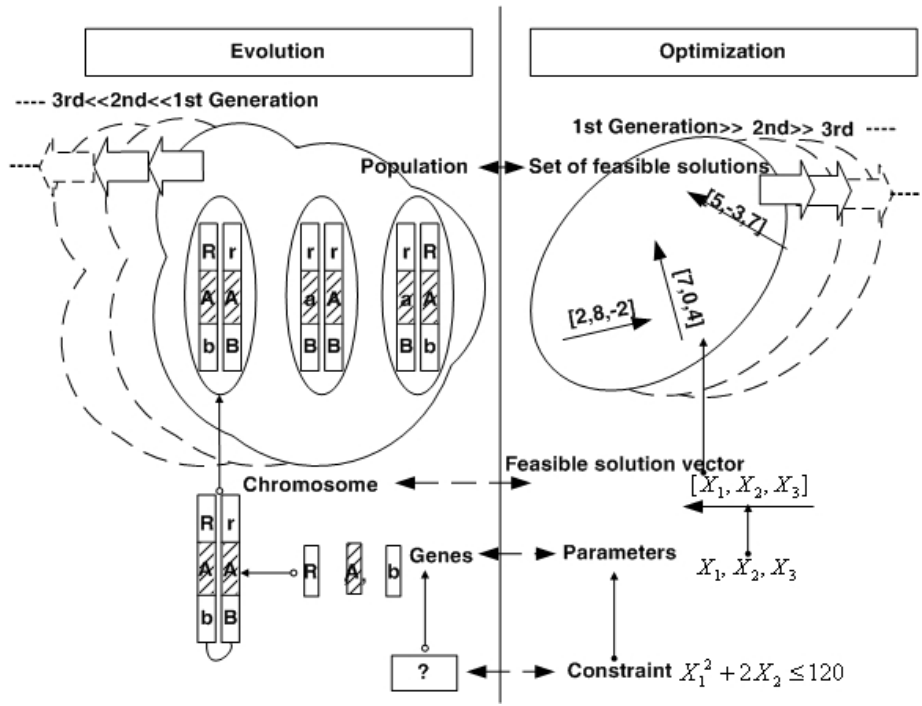

**Fig. 2.** Analogy between evolution and optimization.

What is the counterpart to constraint in biology?

Satisfaction of a constraint is dependent on the values of the corresponding parameters. Similarly, a group of genes can determine whether a fatal characteristic will be expressed. A characteristic, as its biological definition, is any discriminated phenotypic trait of an organism, generally determined by genes. Hence, we can make a hypothesis:

**A constraint is a selection of characteristics.**

**III.1 Natural Selection**

To push constraints into the objective function, is to affect the selection of solution vectors through making a new "natural law". This is generally known as a penalty method. The algorithm below is to be applied recursively.

**Algorithm 1:** An improved objective function is achieved with a penalty method,

$$\text{Min} \sum_{i=1}^{N}\left[\left(X_d^i - X^i\right)^2 + \left(Y_d^i - Y^i\right)^2\right] + W_1 \cdot P_{\text{f-Gra}} + W_2 \cdot P_{\text{f-Seq}} \qquad (5)$$

In which the penalty functions are

$$P_{\text{f-Gra}} = \begin{cases} 1 & \text{Grashof is not satisfied} \\ 0 & \text{Grashof is satisfied} \end{cases} \quad W_1 = 1000;$$

And

$$P_{\text{f-Seq}} = \begin{cases} 1 & \text{Sequence is not satisfied} \\ 0 & \text{Sequence is satisfied} \end{cases} \quad W_2 = 1000.$$

In Figure 2, equation (3) is substituted by equation (5) as the new objective function in the blocks after R1 and R2. It is named the "natural selective improvement" (NSI) method for constraints.

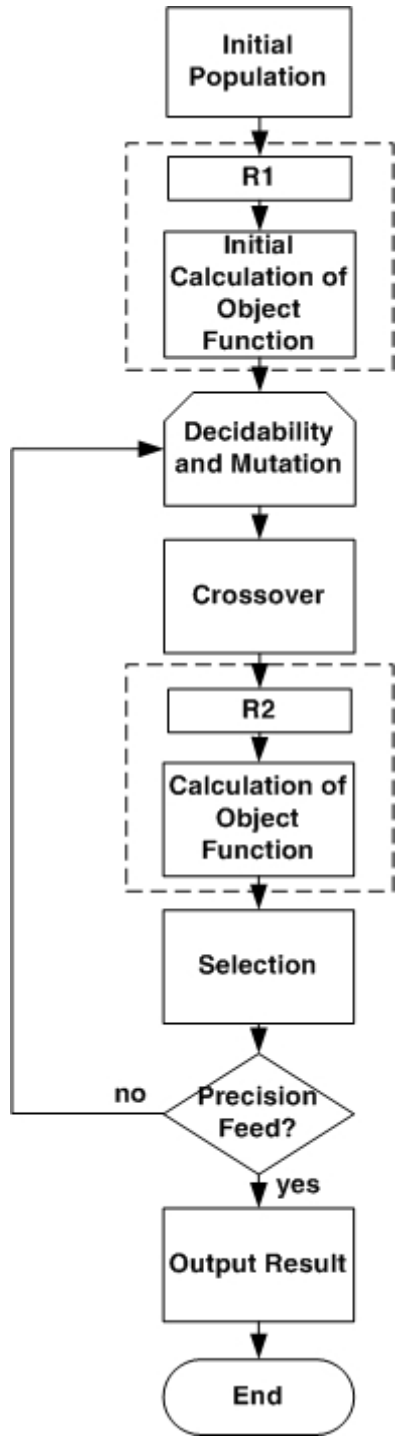

**Fig. 3.** Optimization flow of the differential evolution method. In a typical search, there is no manipulation in blocks R1 and R2. Different manipulations could be put into blocks R1 and R2 and the calculation blocks after them (marked with dashed lines) for better performance.

**III.2 Artificial selection**

Instead of naturally evolving so many species on a scale of millions of years, people achieve their purpose in a faster way. A race horse with sprinting ability, a beautiful goldfish, and a good hunting hound can all be bred quickly in filtered populations. The desired offspring could appear more rapidly if the frequency of evolution in the population is changed.

**Algorithm 2:**

If pre-selection to satisfy a constraint could be accomplished before the beginning of optimization, then a shorter search than NSI can be expected.

In block R1, washing and regenerating a new chromosome can ensure that every randomized chromosome in the initial population is acceptable to the constraints.

**Algorithm 3:**

To make a shorter search, the same work could be executed in every generation. For colleagues in the field of biology, a giant backup population is needed to do the same thing. For us, it may require an indeterminate amount of time. Therefore, a limitation must be placed on update times. NSI will still be reliable if a 100% satisfaction of constraints is impossible in a generation.

These two algorithms are denoted as artificial selective improvement in initial generation (ASI-IG) and artificial selective improvement in all generations (ASI-AG) respectively.

**III.3 Limited selection**

A population would evolve in a limited environment if the states of a gene could be controlled artificially to prevent undesired characteristics. This means that the chromosomes of a population can satisfy the constraints before it is born. Although complete (entirely random) diversity is replaced by a limitation of characteristics, it avoids the passive nature of artificial selection. An invalid chromosome will be substituted by a suitable one immediately.

**Algorithm 4:**

To maintain diversity, a set of pseudo-randomization generators is placed in blocks R1 and R2. These generators can produce "random" chromosomes that satisfy certain constraints. This algorithm is called limited selective improvement (LSI).

**III.4 Lamarckism**

Unlike the natural selection theory of Darwinism, this doctrine, which has been proved incorrect, holds that an individual could actively switch a characteristic by itself in order to adapt. For example, the long neck of a giraffe would not be the result of gradual evolution by longer necks being naturally selected in every generation. According to Lamarckism, an internal driving force made every giraffe try to reach the higher and fresher leaves. These efforts made their necks a little longer. Changes called acquired characteristics would be inherited by the next generation, and after many generations would result in the present-day giraffe.

**Algorithm 5:**

Based on Lamarckism, an algorithm different from ASI and LSI can be designed. Instead of producing a new individual, a modification or local adjustment of an existing unsuitable one is made to meet the constraint condition.

This is called self-selective improvement (SSI).

Because the modification is active and constructive, many kinds of realization are easy to attain. For example, six rules, G1 to G6, have been developed. One of them, G4, is designed as follows,

$$\begin{cases} r_2' = r_2 + S/2 \\ r_3' = r_3 + S/2 + 0.01 \end{cases} \quad (6)$$

Where,

$$S = r_1 - r_4 - r_2 + r_3$$

Then, a modified individual with bars $r_1$, $r_2$', $r_3$', and $r_4$, is constructed to satisfy the Grashof condition.

**IV. Results**

To compare the performance of the presented algorithms, three cases, the same as those of Cabrera's paper [25] (Marked CP in this article) were computed. The last two cases can also be seen in Kunjur and Krishnamurthy's paper [22] (KKP). All results were calculated in MatLab 7.0. Codes were reprogrammed, based on the original codes from webpage [27].

In MatLab environment, five proposed algorithms are coded. For NSI, there is nothing but replacing old objective function with one including penalty items. For ASI-IG and ASI-AG, a piece of circulation program is run with "while" control to ensure an ASI rule. To control time, top times of update to unfeasible chromosomes are 6 in ASI-AG.

Codes of LSI in MatLab are like following as a chromosome generator for Grashof condition.

*function GrashofV =   GrashofP()   %Return Random Grashof vector。*

*%Make a fake random series to ensure satisfaction of Grashof condition*

   *r_max = rand;*

   *r_2    = rand*r_max;*

   *r_3    = r_max+(rand-1)*r_2;*

   *r_min = rand*(r_2+r_3-r_max);*

   *GrashofVpre = [r_max,r_2,r_3,r_min];*

   *GrashofV     = GrashofVpre(randperm(4));*

   *%end*

*Case 1:*

Design a mechanism with chromosome $X = [r_1, r_2, r_3, r_4, r_{ex}, r_{ey}, \theta_0, x_0, y_0, \theta_1^1, ..., \theta_1^6]$ to realize positions $[(20,20), (20,25), (20,30), (20,35), (20,40), (20,45)]$ ,

where $r_1 \sim r_4 \in [0,60]$, $r_{ex}, r_{ey}, x_0, y_0 \in [-60,60]$ and $\theta_0, \theta_1^1 \sim \theta_1^6 \in [0,2\pi]$.

The configuration of DE is: I_NP = 100, I_Itermax =1000, F_weight = 0.3, F_CR = 0.8/1.0, MP = 0.1 and I_strategy = 6.

The best result is $X$ = [40.061, 10.785, 24.47, 43.887, 32.236, 10.064, 3.7921, -2.4468, 56.545, 1.9659, 2.5047, 2.9448, 3.3791, 3.8469, 4.3841] with final error 0.00414.

Calculations were carried out with the five algorithms, and data about them is listed in Table 1. Compare to published results, best convergence of searching and mechanism position are drawn in Figure 4 and Figure 5. A statistic about precision of searching of every algorithm is listed in Table 2.

**Table 1.** Final error and convergence for each algorithm. Errors are recorded when F_CR = 0.8, while convergences in the 100th generation are recorded when F_CR = 1.0.

|  | NSI | ASI-IG | ASI-AG | LSI | SSI | CP |
| --- | --- | --- | --- | --- | --- | --- |
| Final error (mm) | 0.107 | 0.0667 | 0.0276 | 0.00976 | **0.00414[1]** | 0.0363* |
| Convergence (%) | 0 | 98.30 | 83.13 | 99.21 | **99.69[2]** | 99.99 |

Note:

*：Actual error of best result in CP is 0.0363 rather than 0.02617 as reported in the paper.

1：Stopped at 692nd generation.

2：Stopped at 396th generation with acceptable error 0.00627.

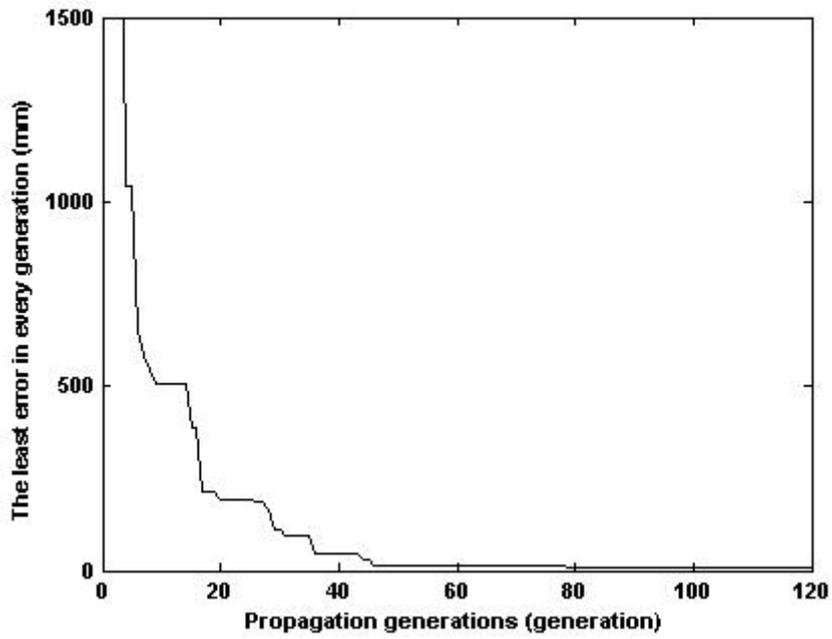

**Fig. 4.** Evolution process of Case 1's objective function. Best convergence was achieved when searched with SSI.

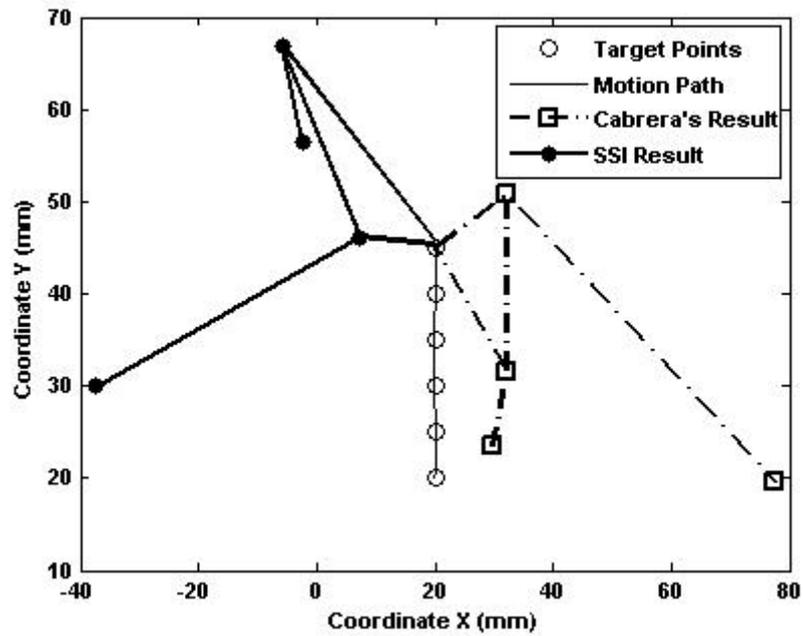

**Fig. 5.** Optimal mechanisms of Case 1. Best mechanism with the least final error is plotted. So did CP's mechanism. Two mechanism have opposite motion directions.

**Table 2.** Statistic of final errors in hundred searching of every algorithm

| Precision (mm) | NSI (times) | ASI-IG (times) | ASI-AG (times) | LSI (times) | SSI (times) |
| --- | --- | --- | --- | --- | --- |
| <0.1 | 1 | 0 | 0 | **5** | 0 |
| <1 | 1 | 5 | 10 | 36 | 13 |
| <10 | 9 | 18 | 45 | 84 | 70 |
| <100 | 33 | 39 | 92 | 100 | 100 |
| <200 | 33 | 41 | 95 | 100 | 100 |
| <500 | 36 | 43 | 98 | 100 | 100 |
| <1000 | 36 | 45 | 99 | 100 | 100 |

*Case 2:*

Solve $X = [r_1, r_2, r_3, r_4, r_{ex}, r_{ey}]$ to pass through the points $[(3.000, 3.000), (2.759, 3.363), (2.372, 3.663), (1.890, 3.862), (1.355, 3.943)]$, when input angle $\theta_1 = [\pi/6, \pi/4, \pi/3, 5\pi/12, \pi/2]$. Where $r_1 \sim r_4 \in [0, 50]$, $r_{ex}, r_{ey} \in [-50, 50]$.

The configuration of DE is: I_NP = 50, I_Itermax =100, F_weight = 0.3, F_CR = 1.0, MP = 0.1 and I_strategy = 6.

The best result is $X$ = [9.3684, 2.0356, 45.497, 43.807, 2.317, -0.27917] with final error 1.354e-4.

Best solutions with every algorithm and consumption time of each search are shown in Table 3. A more restricted range is set to get more precise solutions (see Table 4). A data are gathered in Table 5 to observe the performance of different algorithms. An optimized path is drawn in Figure 6.

Table 3. Best solution of every algorithm. The best result is achieved in calculation with SSI, which is weaker than result of CP's.

|  | NSI | ASI-IG | ASI-AG | LSI | **SSI** | KKP | CP |
|---|---|---|---|---|---|---|---|
| Error (mm) | 5.550e-4 | 4.337e-4 | 7.6215e-4 | 2.283e-4 | **1.354e-4** | 9.526e-4 | 1.827e-6 |
| $r_1$ (mm) | 12.0665 | 11.5679 | 11.9024 | 10.0220 | **9.3684** | 3.5096 | 3.0630 |
| $r_2$ (mm) | 2.0604 | 2.0043 | 2.0912 | 2.0657 | **2.0356** | 1.8576 | 1.9960 |
| $r_3$ (mm) | 37.6264 | 22.7569 | 16.9863 | 21.7393 | **45.4974** | 4.7258 | 3.3058 |
| $r_4$ (mm) | 35.5299 | 18.5983 | 13.4742 | 19.3834 | **43.8072** | 3.5187 | 2.5247 |
| $r_{ex}$ (mm) | 2.2970 | 2.3243 | 2.2377 | 2.3047 | **2.3170** | 1.9595 | 1.6452 |
| $r_{ey}$ (mm) | -0.2753 | 0.4146 | 0.4610 | 0.1101 | **-0.2792** | 1.5589 | 1.7090 |
| Time(s) | 1.11 | 1.05 | 1.25 | 1.20 | **1.23** | 16.98 | 2.86 |

Table 4. More precise results for Case 2. Best solution of every algorithm is achieved within parameter range [0, 5]. Calculation with each algorithm is fulfilled without else modification. Under this situation, all results caused by calculation with presented algorithms are better than one of CP.

|  | NSI | ASI-IG | ASI-AG | **LSI** | SSI | **CP** |
|---|---|---|---|---|---|---|
| Error (mm) | 7.7785e-007 | 7.5832e-007 | 7.5996e-007 | **7.5631e-007** | 7.5711e-007 | **1.8270e-006** |
| $r_1$ (mm) | 3.6770 | 3.6622 | 3.6818 | **3.6470** | 3.6758 | **3.0630** |
| $r_2$ (mm) | 1.9971 | 1.9979 | 1.9984 | **1.9978** | 1.9980 | **1.9960** |
| $r_3$ (mm) | 4.0077 | 3.9866 | 4.0079 | **3.9699** | 4.0031 | **3.3058** |
| $r_4$ (mm) | 2.7097 | 2.7012 | 2.7067 | **2.6976** | 2.7072 | **2.5247** |
| $r_{ex}$ (mm) | 1.6736 | 1.6713 | 1.6716 | **1.6710** | 1.6724 | **1.6452** |
| $r_{ey}$ (mm) | 1.6796 | 1.6807 | 1.6798 | **1.6812** | 1.6796 | **1.7090** |

**Table 5.** Statistic of final errors in hundred searching of every algorithm.

| Precision (mm) | NSI (times) | ASI-IG (times) | ASI-AG (times) | LSI (times) | SSI (times) |
|---|---|---|---|---|---|
| $<10^{-5}$ | 0 | 0 | 0 | 0 | 0 |
| $<10^{-4}$ | 0 | 0 | 0 | 0 | 0 |
| $<10^{-3}$ | 2 | 2 | 1 | **11** | 3 |
| $<10^{-2}$ | 57 | 55 | 69 | 78 | 53 |
| $<10^{-1}$ | 57 | 56 | 70 | 79 | 53 |
| $<10^{0}$ | 97 | 87 | 94 | 100 | 90 |
| $<10^{1}$ | 100 | 100 | 100 | 100 | 100 |

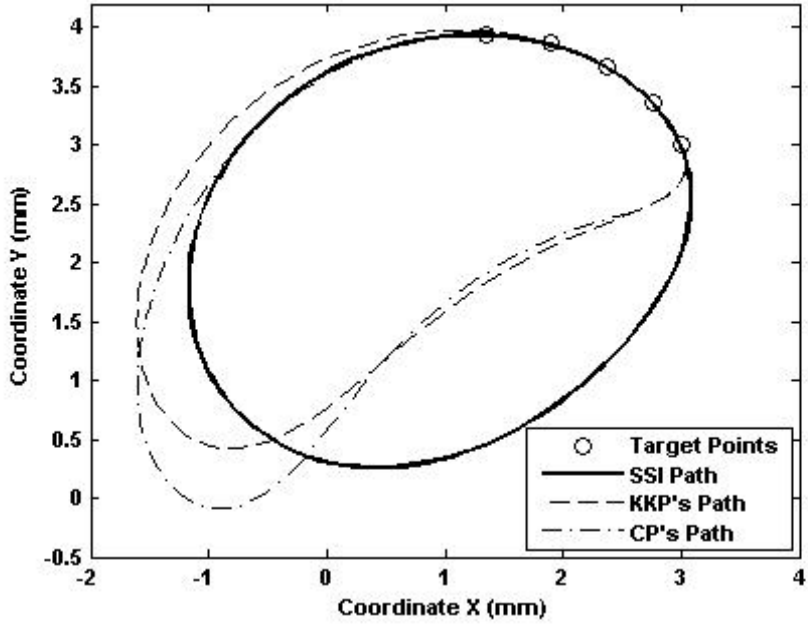

**Fig. 6.** Optimal paths of Case 2. Best path with the least final error is plotted. So did CP's and KKP's mechanism. An amazing result was reported in CP.

*Case 3:*

Determine $X = [r_1, r_2, r_3, r_4, r_{ex}, r_{ey}, \theta_0, x_0, y_0, \theta_1^1]$ to pass through the points [(0.5, 1.1), (0.4, 1.1), (0.3, 1.1), (0.2, 1.0), (0.1, 0.9), (0.005, 0.75), (0.02, 0.6), (0.0, 0.5), (0.0, 0.4), (0.03, 0.3), (0.1, 0.25), (0.15, 0.2), (0.2, 0.3), (0.3, 0.4), (0.4, 0.5), (0.5, 0.7), (0.6, 0.9), (0.6, 1.0)], when input angle $\theta_1 = [\theta_1^1, \theta_1^1 + \dfrac{\pi}{9} \times i]$ (*i* = 1~17), where $r_1 \sim r_4 \in [0, 50]$, $r_{ex}, r_{ey}, x_0, y_0 \in [-50, 50]$.

Configuration of DE is: I_NP = 100, I_Itermax = 50, F_weight = 0.05, F_CR = 1.0, MP = 0.1 and I_strategy = 6.

The best result is *X* = [25.97, 0.46768, 33.636, 34.753, 3.7388, -13.005, 4.0098, 9.6911, 10.36, 3.0223] with final error 0.0476.

Best solutions with every algorithm and consumption time of each search are shown in Table 6. Also, a path compare is shown in Figure 7. A statistic table is achieved. (See Table 7)

**Table 6.** Final error and time of calculation for case 3.

|  | NSI | ASI-IG | ASI-AG | LSI | SSI | KKP | CP |
|---|---|---|---|---|---|---|---|
| Error (mm) | 0.0514 | 0.0863 | **0.0476** | 0.0498 | 0.0568 | 0.0653* | 0.0348** |
| $r_1$ (mm) | 11.0968 | 27.4608 | **25.9700** | 10.8045 | 43.7593 | 1.8797 | 3.0579 |
| $r_2$ (mm) | 0.2604 | 0.2725 | **0.4677** | 0.4435 | 0.2710 | 0.2749 | 0.2378 |
| $r_3$ (mm) | 17.4402 | 30.6331 | **33.6361** | 12.7383 | 45.3331 | 1.1803 | 4.8290 |
| $r_4$ (mm) | 13.9683 | 35.0674 | **34.7532** | 14.5478 | 17.3158 | 2.1382 | 2.0565 |
| $r_{ex}$ (mm) | -10.1608 | 1.2751 | **3.7388** | 0.9218 | -0.6060 | -0.8336 | 0.7670 |
| $r_{ey}$ (mm) | 1.7121 | 16.3477 | **-13.0051** | -4.8433 | 47.5827 | -0.3788 | 1.8508 |
| $\theta_0$ (rad) | 2.2385 | 2.5740 | **4.0098** | 3.9977 | 0.3226 | 4.3542 | 1.0022 |
| $x_0$ (mm) | -9.9504 | -9.6239 | **9.6911** | 3.6881 | 31.8072 | 1.1321 | 1.7768 |
| $y_0$ (mm) | 2.1289 | 13.7719 | **10.3599** | 4.1837 | -34.9621 | 0.6634 | -0.6420 |
| $\theta_1^1$ (rad) | 4.7470 | 4.8284 | **3.0223** | 3.0096 | 1.0261 | 2.5586 | 0.2262 |
| Time(s) | 3.06 | 3.33 | **3.59** | 3.59 | 3.47 | 37.03 | 3.25 |

Note:

*: Actual error of best result in KKP is 0.0653 rather than 0.0430.

**: Actual error of best result in CP is 0.0348 rather than 0.0245.

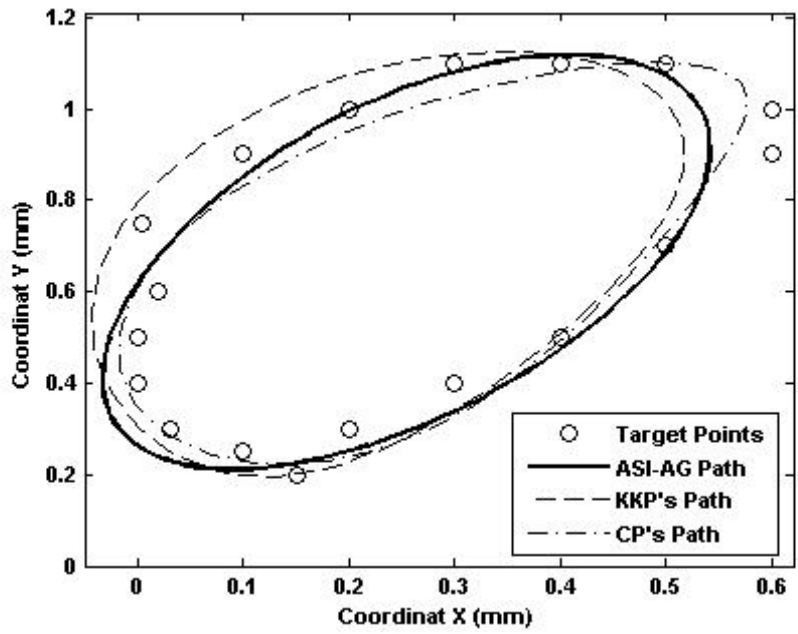

**Fig. 7.** Optimal paths of Case 3. Best path with the least final error is plotted. So did CP's and KKP's mechanism.

**Table 7.** Statistic of final errors in hundred searching of every algorithm

| Precision (mm) | NSI (times) | ASI-IG (times) | ASI-AG (times) | LSI (times) | SSI (times) |
| --- | --- | --- | --- | --- | --- |
| <0.1 | 0 | 0 | 0 | 0 | 0 |
| <1 | 0 | 0 | 0 | **7** | 1 |
| <10 | 2 | 0 | 47 | **92** | 78 |
| <100 | 18 | 17 | 90 | 100 | 90 |
| <200 | 29 | 24 | 94 | 100 | 92 |
| <500 | 52 | 44 | 99 | 100 | 92 |
| <1000 | 65 | 62 | 100 | 100 | 93 |
| =1000 | 35 | 38 | 0 | 0 | 0 |

## V. Discussion

According to Table 1, the convergence of the method in CP is better than that of SSI (See Figure 4), the best performing of all five algorithms. However, a good performance at the 100th generation does not necessarily result in a better consequence. The search with NSI shown in Table 1 recorded a 0 convergence, still achieved a result with final error 1.764. All final errors except NSIs are better than CPs. In addition, a bigger initial error could cause a better 100th convergence.

In case 1, SSI is the only algorithm that stopped before the limitation of 1000 generations. Our best result is achieved at the 692nd generation. The convergence test of SSI stopped at the 396th generation with a final error 0.00627. Its solution is [46.569, 10.1, 26.69, 47.194, 32.356, 12.288, 0.60687, 40.963, 5.8553, 4.6127, 3.8929, 3.36, 2.8914, 2.4015, 1.6726].

For case 2, all final errors of the presented algorithms are better than KKPs, but can not compare with CPs. Actually, it seems a huge challenge to reach that precision in so broad a range. Another calculation (See Table 4) indicates that all five presented algorithms achieved results comparable to CPs when the parameter range was shortened to [0, 5].

Final errors shown in Table 3 and Table 4 demonstrate close performance of these algorithms.

According to Figure 2, the last four algorithms deal with constraints before calculation of the objective function. The results support these improvements.

Although development of hardware plays an important role, the calculation speed has increased relative to earlier work.

Among the presented algorithms, LSI and SSI are better than the others in optimization. This is reasonable, according to the ideas behind them.

Generally, the configuration of DE influences the search process. For instance, better results are achieved with F_CR = 0.8 than with F_CR = 1.0 in case 1. On the contrary, convergence of the algorithms is better with F_CR = 1.0 than with F_CR = 0.8.

Hundreds or thousands of searches with each algorithm in each case indicated that LSI and SSI always perform well. LSI always possesses the highest percentage of optimal solutions. SSI often achieves the best solution. (See Table 2, Table 5, and Table 7)

A further study might focus on SSI. It is obvious that active modification via Lamarckism can only build a subset of the feasible set. This will become an obstacle to obtaining the optimal results, although a restricted parameter space is good for searching. Application of the presented algorithms to other optimization problems is necessary to test their performance.

**VI. Conclusion**

One algorithm is re-explained and others are presented to handle constraints in the synthesis of a four-bar mechanism. All of the algorithms originate from evolutionary theory. The algorithms and their methodology are useful for nonlinear constrained optimization problems in all fields.

Calculations of three cases showed better or comparable performance compared with previous methods, together with high speed.

These biologically selective improvements (especially with respect to LSI and SSI) for handling constraints in nonlinear optimization are worthy of further study. Especially, optimization with LSI gives an new and effective means, pseudo-randomization generator, for constraint handling. And SSI presents a most simple way to overcome constraint.

**Table 1.** Final error and convergence for each algorithm. Errors are recorded when F_CR = 0.8, while convergences in the 100th generation are recorded when F_CR = 1.0.

|  | NSI | ASI-IG | ASI-AG | LSI | SSI | CP |
|---|---|---|---|---|---|---|
| Final error (mm) | 0.107 | 0.0667 | 0.0276 | 0.00976 | **0.00414[1]** | 0.0363* |
| Convergence (%) | 0 | 98.30 | 83.13 | 99.21 | **99.69[2]** | 99.99 |

Note:

\*：Actual error of best result in CP is 0.0363 rather than 0.02617 as reported in the paper.

1：Stopped at 692nd generation.

2：Stopped at 396th generation with acceptable error 0.00627.

**Table 2.** Statistic of final errors in hundred searching of every algorithm

| Precision (mm) | NSI (times) | ASI-IG (times) | ASI-AG (times) | LSI (times) | SSI (times) |
| --- | --- | --- | --- | --- | --- |
| <0.1 | 1 | 0 | 0 | **5** | 0 |
| <1 | 1 | 5 | 10 | 36 | 13 |
| <10 | 9 | 18 | 45 | 84 | 70 |
| <100 | 33 | 39 | 92 | 100 | 100 |
| <200 | 33 | 41 | 95 | 100 | 100 |
| <500 | 36 | 43 | 98 | 100 | 100 |
| <1000 | 36 | 45 | 99 | 100 | 100 |

Table 3. Best solution of every algorithm. The best result is achieved in calculation with SSI, which is weaker than result of CP's.

|  | NSI | ASI-IG | ASI-AG | LSI | **SSI** | KKP | CP |
|---|---|---|---|---|---|---|---|
| Error (mm) | 5.550e-4 | 4.337e-4 | 7.6215e-4 | 2.283e-4 | **1.354e-4** | 9.526e-4 | 1.827e-6 |
| $r_1$ (mm) | 12.0665 | 11.5679 | 11.9024 | 10.0220 | **9.3684** | 3.5096 | 3.0630 |
| $r_2$ (mm) | 2.0604 | 2.0043 | 2.0912 | 2.0657 | **2.0356** | 1.8576 | 1.9960 |
| $r_3$ (mm) | 37.6264 | 22.7569 | 16.9863 | 21.7393 | **45.4974** | 4.7258 | 3.3058 |
| $r_4$ (mm) | 35.5299 | 18.5983 | 13.4742 | 19.3834 | **43.8072** | 3.5187 | 2.5247 |
| $r_{ex}$ (mm) | 2.2970 | 2.3243 | 2.2377 | 2.3047 | **2.3170** | 1.9595 | 1.6452 |
| $r_{ey}$ (mm) | -0.2753 | 0.4146 | 0.4610 | 0.1101 | **-0.2792** | 1.5589 | 1.7090 |
| Time(s) | 1.11 | 1.05 | 1.25 | 1.20 | **1.23** | 16.98 | 2.86 |

**Table 4.** More precise results for Case 2. Best solution of every algorithm is achieved within parameter range [0, 5]. Calculation with each algorithm is fulfilled without else modification. Under this situation, all results caused by calculation with presented algorithms are better than one of CP.

|  | NSI | ASI-IG | ASI-AG | LSI | SSI | **CP** |
|---|---|---|---|---|---|---|
| Error (mm) | 7.7785e-007 | 7.5832e-007 | 7.5996e-007 | **7.5631e-007** | 7.5711e-007 | **1.8270e-006** |
| $r_1$ (mm) | 3.6770 | 3.6622 | 3.6818 | **3.6470** | 3.6758 | **3.0630** |
| $r_2$ (mm) | 1.9971 | 1.9979 | 1.9984 | **1.9978** | 1.9980 | **1.9960** |
| $r_3$ (mm) | 4.0077 | 3.9866 | 4.0079 | **3.9699** | 4.0031 | **3.3058** |
| $r_4$ (mm) | 2.7097 | 2.7012 | 2.7067 | **2.6976** | 2.7072 | **2.5247** |
| $r_{ex}$ (mm) | 1.6736 | 1.6713 | 1.6716 | **1.6710** | 1.6724 | **1.6452** |
| $r_{ey}$ (mm) | 1.6796 | 1.6807 | 1.6798 | **1.6812** | 1.6796 | **1.7090** |

**Table 5.** Statistic of final errors in hundred searching of every algorithm.

| Precision (mm) | NSI (times) | ASI-IG (times) | ASI-AG (times) | LSI (times) | SSI (times) |
|---|---|---|---|---|---|
| $<10^{-5}$ | 0 | 0 | 0 | 0 | 0 |
| $<10^{-4}$ | 0 | 0 | 0 | 0 | 0 |
| $<10^{-3}$ | 2 | 2 | 1 | **11** | 3 |
| $<10^{-2}$ | 57 | 55 | 69 | 78 | 53 |
| $<10^{-1}$ | 57 | 56 | 70 | 79 | 53 |
| $<10^{0}$ | 97 | 87 | 94 | 100 | 90 |
| $<10^{1}$ | 100 | 100 | 100 | 100 | 100 |

**Table 6.** Final error and time of calculation for case 3.

|  | NSI | ASI-IG | ASI-AG | LSI | SSI | KKP | CP |
|---|---|---|---|---|---|---|---|
| Error (mm) | 0.0514 | 0.0863 | **0.0476** | 0.0498 | 0.0568 | 0.0653* | 0.0348** |
| $r_1$ (mm) | 11.0968 | 27.4608 | **25.9700** | 10.8045 | 43.7593 | 1.8797 | 3.0579 |
| $r_2$ (mm) | 0.2604 | 0.2725 | **0.4677** | 0.4435 | 0.2710 | 0.2749 | 0.2378 |
| $r_3$ (mm) | 17.4402 | 30.6331 | **33.6361** | 12.7383 | 45.3331 | 1.1803 | 4.8290 |
| $r_4$ (mm) | 13.9683 | 35.0674 | **34.7532** | 14.5478 | 17.3158 | 2.1382 | 2.0565 |
| $r_{ex}$ (mm) | -10.1608 | 1.2751 | **3.7388** | 0.9218 | -0.6060 | -0.8336 | 0.7670 |
| $r_{ey}$ (mm) | 1.7121 | 16.3477 | **-13.0051** | -4.8433 | 47.5827 | -0.3788 | 1.8508 |
| $\theta_0$ (rad) | 2.2385 | 2.5740 | **4.0098** | 3.9977 | 0.3226 | 4.3542 | 1.0022 |
| $x_0$ (mm) | -9.9504 | -9.6239 | **9.6911** | 3.6881 | 31.8072 | 1.1321 | 1.7768 |
| $y_0$ (mm) | 2.1289 | 13.7719 | **10.3599** | 4.1837 | -34.9621 | 0.6634 | -0.6420 |
| $\theta_1^1$ (rad) | 4.7470 | 4.8284 | **3.0223** | 3.0096 | 1.0261 | 2.5586 | 0.2262 |
| Time(s) | 3.06 | 3.33 | **3.59** | 3.59 | 3.47 | 37.03 | 3.25 |

Note:

*: Actual error of best result in KKP is 0.0653 rather than 0.0430.

**: Actual error of best result in CP is 0.0348 rather than 0.0245.

**Table 7.** Statistic of final errors in hundred searching of every algorithm

| Precision (mm) | NSI (times) | ASI-IG (times) | ASI-AG (times) | LSI (times) | SSI (times) |
| --- | --- | --- | --- | --- | --- |
| <0.1 | 0 | 0 | 0 | 0 | 0 |
| <1 | 0 | 0 | 0 | **7** | 1 |
| <10 | 2 | 0 | 47 | **92** | 78 |
| <100 | 18 | 17 | 90 | 100 | 90 |
| <200 | 29 | 24 | 94 | 100 | 92 |
| <500 | 52 | 44 | 99 | 100 | 92 |
| <1000 | 65 | 62 | 100 | 100 | 93 |
| =1000 | 35 | 38 | 0 | 0 | 0 |

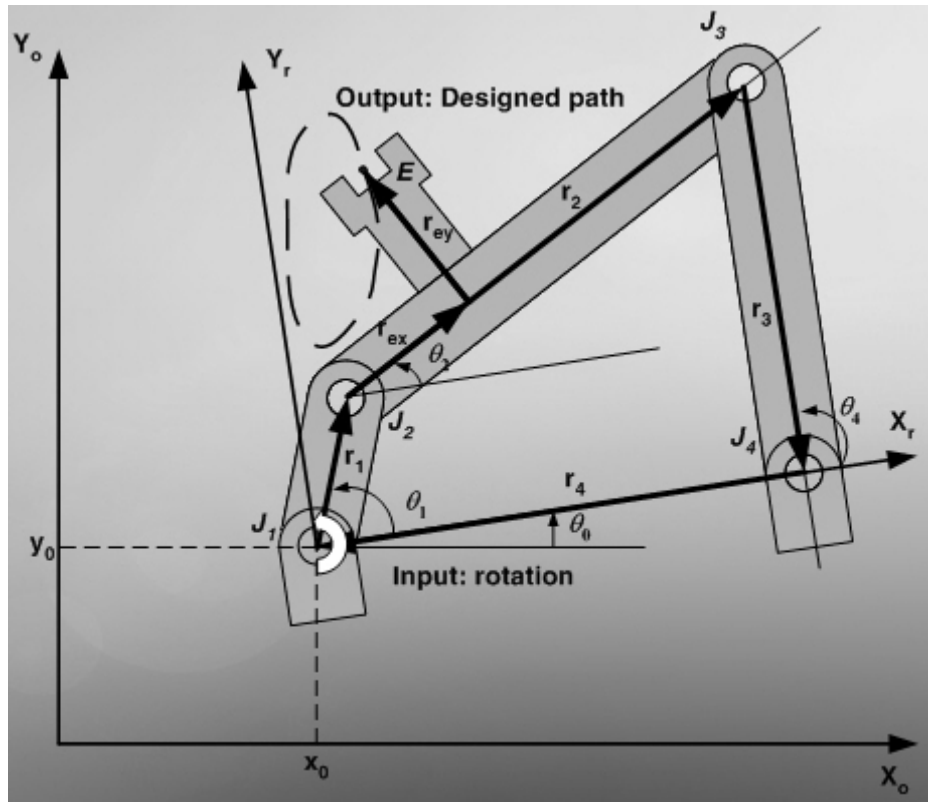

**Fig. 1.** Sketch of a four-bar mechanism. A desired motion path of the endpoint $E$ can be realized from an input rotation in joint $J_1$ via the mechanism. The fourth bar $r_4$ always represents fixed ground.

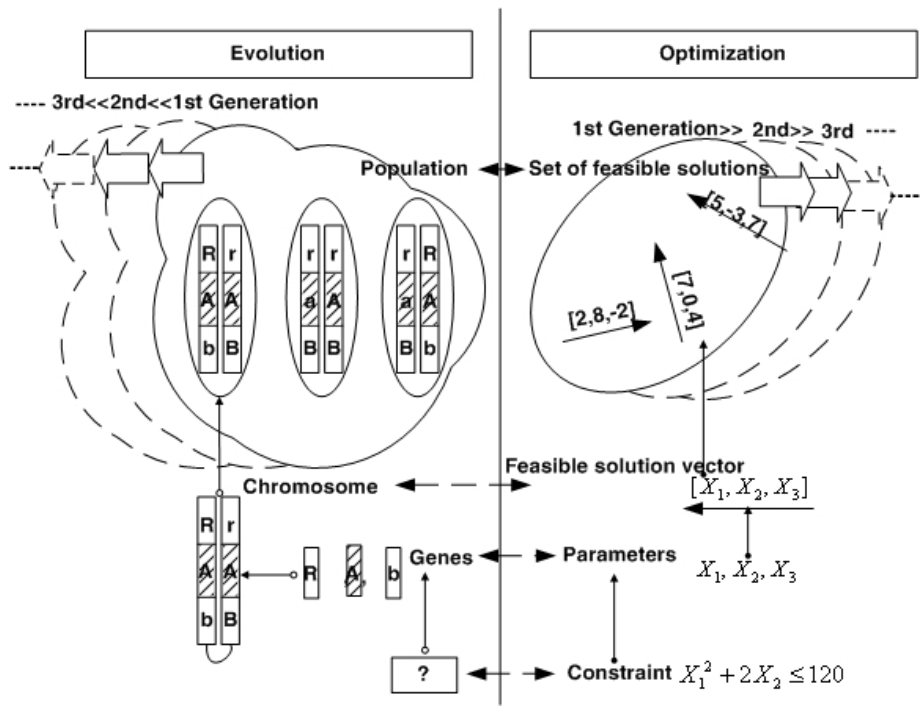

**Fig. 2.** Analogy between evolution and optimization.

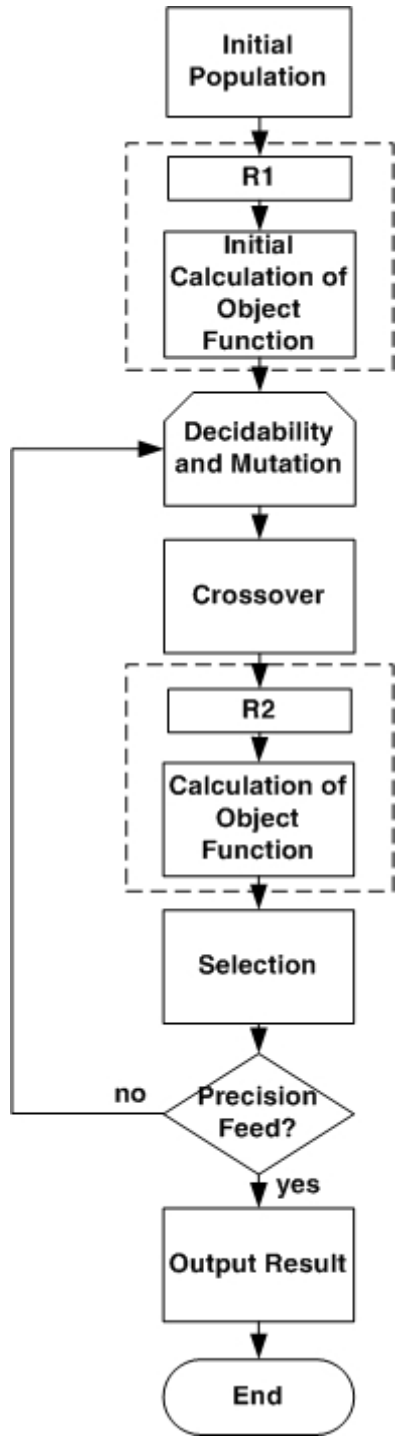

**Fig. 3.** Optimization flow of the differential evolution method. In a typical search, there is no manipulation in blocks R1 and R2. Different manipulations could be put into blocks R1 and R2 and the calculation blocks after them (marked with dashed lines) for better performance.

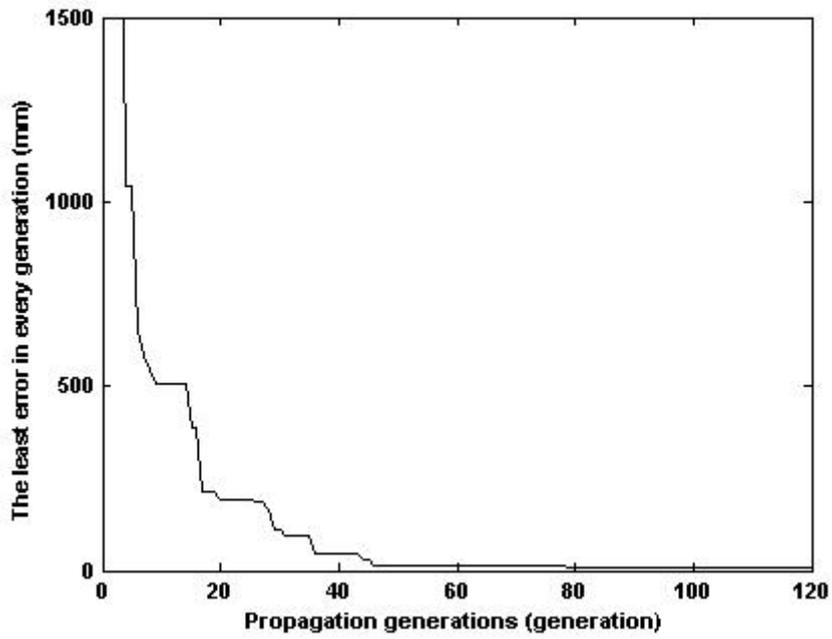

**Fig. 4.** Evolution process of Case 1's objective function. Best convergence was achieved when searched with SSI. The least error of 100[th] generation decreased 99.6912%. It is weaker than CP's result. However, initial least error is a factor.

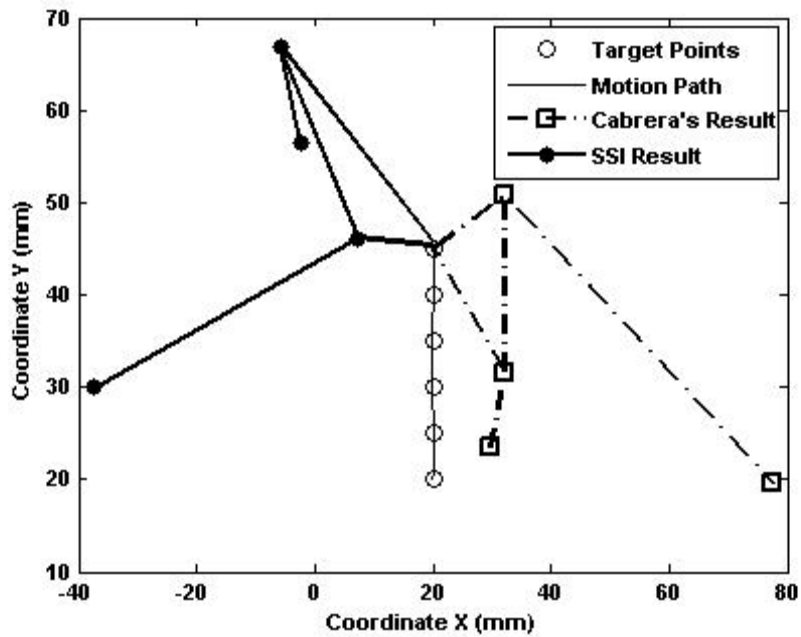

**Fig. 5.** Optimal mechanisms of Case 1. Best mechanism with the least final error is plotted. So did CP's mechanism. Two mechanism have opposite motion directions.

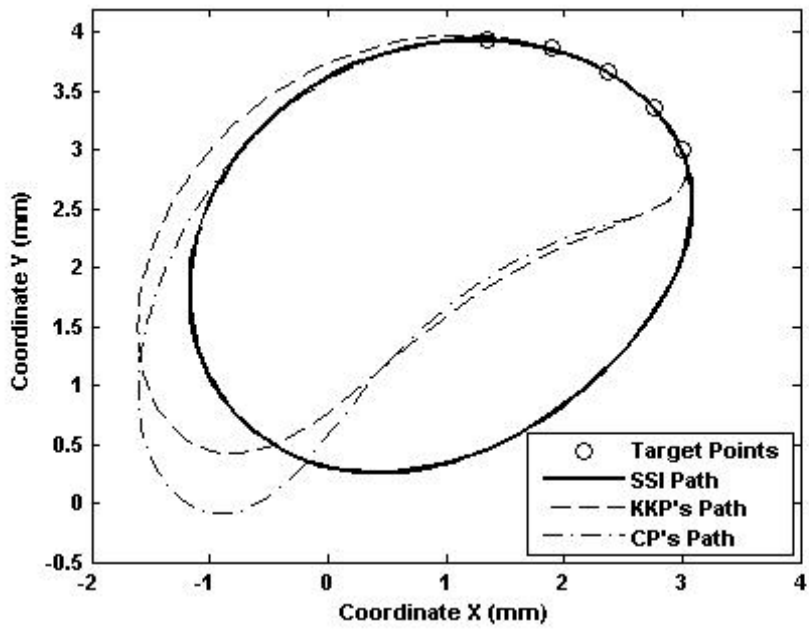

**Fig. 6.** Optimal paths of Case 2. Best path with the least final error is plotted. So did CP's and KKP's mechanism.

An amazing result was reported in CP.

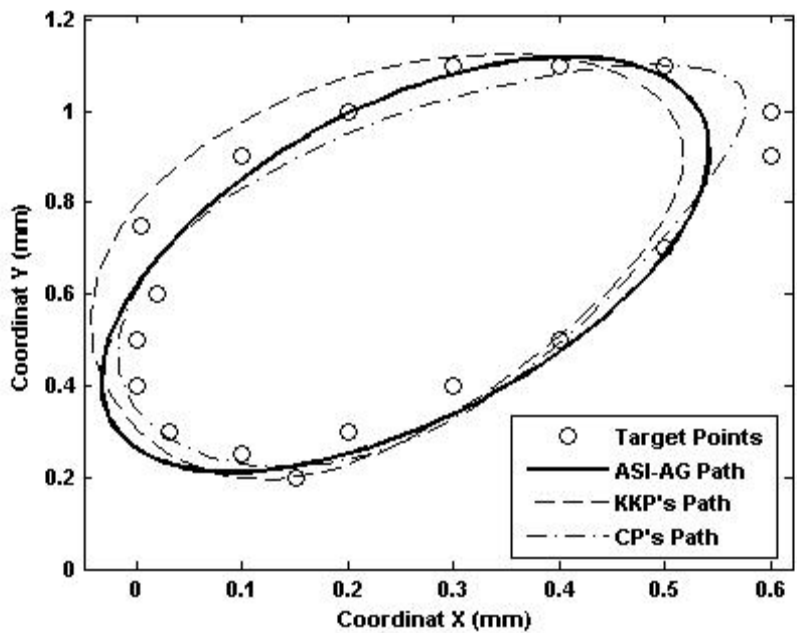

**Fig. 7.** Optimal paths of Case 3. Best path with the least final error is plotted. So did CP's and KKP's mechanism.

## Acknowledgment

This paper is supported by the Initiative Foundation of Scientific Research in CAUC (No. 07qd08s). The author wishes to express his deep appreciation for Dr. Fang Shimin. Without his books popularizing science, the present author would have remained unacquainted with recent ideas in the field of biology. He has changed the opinions of many Chinese with regard to the popularization of science.